\documentclass[a4paper,12pt]{article}
\usepackage{amsmath,amsthm,amsfonts,amssymb,bm,mathrsfs}
\usepackage{enumerate}
\usepackage{graphicx}
\usepackage{natbib}
\usepackage{xcolor}
\usepackage{hyperref}
\usepackage{sgame}

\usepackage{setspace}
\allowdisplaybreaks

\newtheorem{defn}{Definition}
\newtheorem{thm}{Theorem}
\newtheorem{lem}{Lemma}

\newtheorem{prop}{Proposition}
\newtheorem{rmk}{Remark}
\newtheorem{exam}{Example}
\newtheorem{coro}{Corollary}

\newcommand{\bN}{\mathbb{N}}
\newcommand{\bR}{\mathbb{R}}
\newcommand{\cB}{\mathcal{B}}

\newcommand{\cG}{\mathcal{G}}

\newcommand{\cI}{\mathcal{I}}

\newcommand{\cM}{\mathcal{M}}

\newcommand{\cT}{\mathcal{T}}

\newcommand{\ls}{\mbox{Ls}}

\DeclareMathOperator{\rmd}{d\!}

\DeclareMathOperator{\supp}{supp}

\hypersetup{colorlinks=true,
            citecolor=blue,
            linkcolor=blue,
            urlcolor=blue,
            bookmarksopen=true,
            pdfstartview={XYZ null null 1.00},
            pdfpagelayout={SinglePage}
}

\begin{document}
\title{Conditional Expectations  of Correspondences}
\author{Wei~He\thanks{Department of Mathematics, National University of Singapore, 10 Lower Kent Ridge Road, Singapore 119076. E-mail: he.wei2126@gmail.com.}
\and
Yeneng~Sun\thanks{Department of Economics, National University of Singapore, 1 Arts Link, Singapore 117570. Email: ynsun@nus.edu.sg}
}
\date{First version: Nov 19, 2013. This version: \today}
\maketitle

\centerline{Incomplete Preliminary Draft}


\abstract{We characterize the properties of convexity, compactness and preservation of upper hemicontinuity for conditional expectations of correspondences. These results are then applied to obtain a necessary and sufficient condition for the existence of pure strategy equilibria in Bayesian games.}


\newpage
\tableofcontents
\newpage

\section{Introduction}\label{sec-intro}

For a correspondence $F$ from an atomless probability space $(T,\cT,\lambda)$ to the Euclidean space $\bR^n$, let $\cI^{(\cT,\cG)}_{F}$ be the set of all $E(f|\cG)$ such that  $f$ is an integrable selection of $F$, where $E(f|\cG)$ is  the conditional expectation of $f$ on some given sub-$\sigma$-algebra $\cG$ of $\cT$.  By the classical Lyapunov Theorem, the atomless property of $(T,\cT,\lambda)$  implies the convexity of $\cI^{(\cT,\cG)}_{F}$ if $\cG$ is the trivial $\sigma$-algebra $\{T, \emptyset\}$. However, it is easy to see that such a convexity result fails when we work with a general sub-$\sigma$-algebra $\cG$ of $\cT$.\footnote{Suppose that $F(t) = \{0,1\}$ for all $t\in T$, which is a measurable, integral bounded and compact valued correspondence. Given $\cG= \cT$, then the conditional expectation of $F$ conditional on $\cG$ is the set of all integrable selections, which is not convex.} Similarly, some other common regularity properties, such as compactness and preservation of upper hemicontinuity, also fail to hold in the general case. The purpose of this paper is to characterize the properties of convexity, compactness and preservation of upper hemicontinuity for conditional expectations  of correspondences.

The key condition we will work with is that $\cT$ has no $\cG$-atom, which means that $\cT$ does not coincide with $\cG$ when they are restricted on any non-trivial set in $\cT$. Based on this condition, \cite{DE1976} established the equivalence of $\cI^{(\cT,\cG)}_{F}$ and $\cI^{(\cT,\cG)}_{\mbox{co} (F)}$ for any measurable, integrably bounded and closed valued correspondence $F$, where $\mbox{co}(F) (t)$ is the convex hull of $F(t)$ for each $t \in T$. We show that  $\cI^{(\cT,\cG)}_{F}$ is convex for any correspondence $F$ if and only if $\cT$ has no $\cG$-atom. We also prove that this condition is necessary and sufficient for the  weak/weak$^*$ compactness of $\cI^{(\cT,\cG)}_{F}$ for any integrably bounded and closed valued correspondence $F$. A similar necessity and sufficiency result holds for the property on preservation of upper hemicontinuity. Thus, we not only generalize the classical results on integration of correspondences\footnote{See, for example, \cite{Aumann1965} and \cite{Hd1974}.} to the case of conditional expectation, but also demonstrate the optimality of the relevant condition.

To illustrate the application of the main results, Bayesian games with finite actions will be considered.\footnote{There is a substantial literature studying the existence of pure strategy equilibrium  in Bayesian games with finite actions, see \cite{RR1982}, \cite{MW1985}, \cite{KRS2006} and \cite{BD2013}.} We formulate the notion of ``inter-player information'' to describe the influence of player~$i$'s private information in other players' payoffs. The condition of ``coarser inter-player information'' is proposed below and we show that this condition is not only sufficient but also necessary for the existence of pure strategy equilibrium. In particular, interdependent payoffs and correlated types are allowed in our setting. In addition, we prove the  purification results for any behavioral strategy profile.

The rest of the paper is organized as follows. Some basic definitions are given in  Section~\ref{sec-basics}. We state the results on conditional expectations of correspondences in Section \ref{sec-results}. An application to Bayesian games is presented in Section~\ref{sec-Bayesian game}. The proofs are left in Section~\ref{sec-proof}.

\section{Basic Definitions}\label{sec-basics}

Suppose that $(T,\cT,\lambda)$ is a complete probability space endowed with a countably additive probability measure $\lambda$. A correspondence $F$ from $T$  to $\bR^n$ is a mapping from $T$ to the family of nonempty subsets of $\bR^n$.  It is said to be measurable if for all open sets $O\subseteq \bR^n$, we have $\{t\in T\colon F(t)\cap O\neq \emptyset\} \in \cT$.
A measurable function $f$ from $T$ to $\bR^n$ is called a selection of $F$ if $f(t) \in F(t)$ for $\lambda$-almost all $t \in T$. A correspondence $F$ from $T$ to $\bR^n$ is said to be convex (resp. closed, compact) valued if $F(t)$ is convex (resp. closed, compact) for $\lambda$-almost all $t\in T$.

A correspondence $F$ from a topological space $Y$ to another  topological space $Z$ is said to be upper hemicontinuous at $y_0\in Y$ if for any open set $O_Z$ that contains $F(y_0)$, there exists an open neighborhood $O_Y$ of $y_0$ such that $\forall y\in O_Y$, $F(y)\subseteq O_Z$.
$F$ is upper hemicontinuous if it is upper hemicontinuous at every point $y\in Y$.

Suppose that $\cG$ is a sub-$\sigma$-algebra of $\cT$.
For any correspondence $F$ from $T$ to $\bR^n$, let
$$\cI^{(\cT,\cG)}_{F}=\{E(f|\cG)\colon f\mbox{ is an integrable selection of } F  \},
$$
where the conditional expectation  is taken with respect to the probability measure $\lambda$.

Let $L^{\cG}_p(T,\bR^n)$ and $L^{\cG}_{\infty}(T,\bR^n)$  be the set of all $\cG$-measurable mappings from $T$ to $\bR^n$ with the usual norm. That is,
$$L^{\cG}_p(T,\bR^n) = \left\{f\colon f \mbox{ is } \cG \mbox{-measurable and } \left( \int_T \|f\|^p \rmd \lambda  \right)^{\frac{1}{p}} < \infty \right\},
$$
$$L^{\cG}_{\infty}(T, \bR^n) = \{f\colon f \mbox{ is } \cG \mbox{-measurable and essentially bounded under } \lambda \},
$$
where $\|\cdot\|$ is the usual norm in $\bR^n$ and $1\le p <\infty$.  By the Riesz Representation Theorem (see Theorem~13.26-28 of \cite{AB2006}), $L^{\cG}_{q}(T, \bR^n)$ can be viewed as
the dual space of  $L_{p}^\cG(T,\bR^n)$, where $\frac{1}{p} + \frac{1}{q} = 1$. Similarly, one can define $L^{\cT}_p(T,\bR^n)$ and $L^{\cT}_{\infty}(T,\bR^n)$.

If there is a real valued function $h \in L^{\cT}_p(T,\bR)$  such that $\sup\{\|x\|: x\in F(t)\} \leq h(t)$ for $\lambda$-almost all $t\in T$, then the correspondence  $F$ is said to be $p$-integrably bounded, $1\le p \le \infty$ (it is said to be integrably bounded if $p = 1$).

For any nonnegligible subset $D\in\cT$, the restricted probability space $(D,\cG^D,\lambda^D)$ is defined as follows: $\cG^D$ is the $\sigma$-algebra $\{ D\cap D' \colon D'\in\cG \}$ and $\lambda^D$ the probability measure re-scaled from the restriction of $\lambda$ on $\cG^D$. Furthermore, $(D,\cT^D,\lambda^D)$ can be defined similarly. A subset $D\in\cT$ is said to be a $\cG$-atom if $\lambda(D)>0$ and given any $D_0 \in \cT^D$, there exists a set $D_1\in\cG^D$ such that $\lambda(D_0\triangle D_1) = 0$.

\section{The Main Results}\label{sec-results}

In this section,  we will show that the condition that $\cT$ has no $\cG$-atom is sufficient and necessary for the validity of several regularity properties for conditional expectations of correspondences (convexity, compactness and upper hemicontinuity).

The sufficiency part of the following theorem is due to \citet[Theorem~1.2]{DE1976}, while the necessity part is from \citet[Proposition~1]{HS2013}.

\begin{thm}\label{thm-convex hull}
$\cI^{(\cT,\cG)}_F = \cI^{(\cT,\cG)}_{\mbox{co} (F)}$ for any $\cT$-measurable, integrably bounded and closed valued correspondence $F$  if and only if $\cT$ has no $\cG$-atom.
\end{thm}

By the Kuratowski-Ryll-Nardzewski Selection Theorem (see \citet[Theorem~18.13]{AB2006}),  $\cI^{(\cT,\cG)}_F$ is nonempty  for any $\cT$-measurable, integrably bounded and closed valued correspondence $F$.

If $\cT$ has no $\cG$-atom, then the convexity of $\cI^{(\cT,\cG)}_F$ for any correspondence $F$ is a simple corollary of Theorem~\ref{thm-convex hull}. It can be also shown that this condition is necessary for such convexity property.

\begin{coro}\label{coro-convexity}
The set $\cI^{(\cT,\cG)}_F$ is convex for any correspondence $F$ if and only if $\cT$ has no $\cG$-atom.
\end{coro}

Next, we consider the weak/weak$^*$ compactness of $\cI^{(\cT,\cG)}_F$ for a correspondence $F$.

\begin{thm}\label{thm-compact}
The set $\cI^{(\cT,\cG)}_F$ is weakly compact (resp. weak$^*$  compact) in $L_{p}^\cG(T,\bR^n)$ when $1 \le p < \infty$ (resp. $p = \infty$, and $\cG$ is countably generated) for any $p$-integrably bounded and closed valued correspondence $F$ if and only if $\cT$ has no $\cG$-atom.\footnote{The equivalence of compactness and sequential compactness in the weak topology of  a Banach space is important in the proofs of Theorems~\ref{thm-compact} and \ref{thm-upper} below. Such an equivalence still holds for the weak$^*$ topology of $L_{\infty}^\cG(T,\bR^n)$ when $\cG$ is countably generated.}
\end{thm}

The last property is the preservation of weak/weak$^*$ upper hemicontinuity for conditional expectations of correspondences.

\begin{thm}\label{thm-upper}
The following conditions are equivalent.
\begin{enumerate}
  \item For any closed valued correspondence $F$ from $T\times Y\to \bR^n$ ($Y$ is a metric space) such that there is a $p$-integrably bounded and compact valued correspondence $G$ from $T$ to $\bR^n$ and
      \begin{description}
        \item[a.] $F(t,y)\subseteq G(t)$ for $\lambda$-almost all $t\in T$ and all $y\in Y$;
        \item[b.] $F(\cdot,y)$ is $\cT$-measurable for all $y\in Y$;
        \item[c.] $F(t,\cdot)$ is upper hemicontinuous for $\lambda$-almost all $t\in T$;
      \end{description}
      $H(y)= \cI_{F_y}^{(\cT,\cG)}$ is weakly (resp. weak$^*$) upper hemicontinuous in $L_{p}^\cG(T,\bR^n)$ when $1 \le p < \infty$ (resp. $p = \infty$, and $\cG$ is countably generated).
  \item $\cT$ has no $\cG$-atom.
\end{enumerate}
\end{thm}

Note that if $1 < p < \infty$, then $L_p^\cG(T,\bR^n)$ is reflexive. Thus, the weak compactness (resp. weak upper hemicontinuity) and the weak$^*$ compactness (resp. weak$^*$ upper hemicontinuity) are equivalent in $L_p^\cG(T,\bR^n)$ for $1 < p < \infty$.

\begin{rmk}\label{rmk-stochastic games}
\cite{HS2013} proved the existence of stationary Markov perfect equilibria in discounted stochastic games with coarser transition kernels by using Theorem~1.2 of \cite{DE1976}. Recall that  $P$ is an equilibrium payoff correspondence from $T\times V$ to $\bR^n$ such that $P(t,\cdot)$ is upper hemicontinuous and $P(\cdot,v)$ is $\cT$-measurable. Let $R(v)$ be
the set of all selections of $P_v$ for each $v\in V$. The classical Fan-Glicksberg Fixed Point Theorem is applied to the correspondence $\mbox{co}(R)$, which is convex valued and upper hemicontinuous, to obtain a selection $v'$ of $\mbox{co} P(\cdot,v')$. Theorem~\ref{thm-convex hull} then implies the existence of a
selection $v^*$ of $P(\cdot,v')$ such that $E(v'|\cG) = E(v^*|\cG)$, which leads to the existence of a stationary Markov perfect equilibrium. By Corollary~\ref{coro-convexity}, and Theorems~\ref{thm-compact} and \ref{thm-upper}, $\cI_{P_v}^{(\cT,\cG)}$ is  convex, compact valued and upper hemicontinuous. Then the existence result can be also proved by applying the Fan-Glicksberg Fixed Point Theorem to the correspondence $\cI_{P_v}^{(\cT,\cG)}$.
\end{rmk}

\section{Bayesian Games with Inter-player Information}\label{sec-Bayesian game}

In this section, we shall propose the condition of ``coarser inter-player information'', and show that this condition is not only sufficient but also necessary for the existence of pure strategy equilibria in Bayesian games with finite actions. Purification results of behavioral strategy profiles will be also considered.

\subsection{Model}\label{subsec-Bayesian model}

A Bayesian game $\Gamma$  can be described as follows:
\begin{itemize}
  \item The set of players: $I=\{1,2,\ldots,n\}$.
  \item The (private) information space for each player: $\{T_i\}_{i\in I}$. Each $T_i$ is endowed with a countably generated $\sigma$-algebra $\cT_i$. Let $T=\times_{i=1}^nT_i$ and $\cT=\otimes_{i=1}^n\cT_i$.
  \item For each player $i\in I$, $X_i$ is a finite set of actions. Let $X= \prod_{1\leq i \leq n}X_i$.
  \item The information structure: $\lambda$, a probability measure on the measurable space $(T,\cT)$. For each $i \in I$, $\lambda_i$ is the marginal probability of $\lambda$ on $T_i$ and $(T_i,\cT_i,\lambda_i)$ is atomless. $\lambda$ is absolutely continuous with respect to $\bigotimes_{1\le i \le n} \lambda_i$ and $q(t_1,\ldots, t_n)$ is  the Radon-Nikodym derivative.\footnote{This assumption is standard in the literature, see \cite{MW1985}.}
  \item The payoff functions: $\{u_i\}_{i\in I}$. Each $u_i$ is an integrably bounded mapping from $X\times T$ to $\bR$ such that $u_i(x,\cdot)$ is $\cT$-measurable for each $x\in X$.
\end{itemize}

Hereafter,  the notation $-i$ denotes the set of all players except player~$i$. Let $\lambda_{-i} = \otimes_{j\neq i}\lambda_j$. Without loss of generality, we can assume that the mixture of actions of player~$i$ is the simplex $\cM(X_i)$, and the pure actions in $X_i$ correspond to vertices of $\cM(X_i)$. For each player~$i\in I$, a behavioral strategy (resp. pure strategy) is a measurable function from $T_i$ to $\cM(X_i)$ (resp. $X_i$), and $L_i^{\cT_i}$ is the set of all behavioral strategies. $L^{\cT} = \times_{i\in I} L_i^{\cT_i}$.

Given a strategy profile $f=(f_1,\ldots,f_n)$, player~$i$'s expected payoff is
$$U_i(f) = \int_{T} \int_{X} u_i(x, t) \prod_{j\in I} f_j(t_j,\rmd x_j) \lambda (\rmd t).$$

A behavioral (resp. pure) strategy equilibrium is a behavioral (resp. pure) strategy profile $f^*=(f^*_1,f_2^*,\ldots,f^*_n)$ such that $f^*_i$ maximizes $U_i(f_i,f^*_{-i})$  for each $i\in I$.

Consider the \textbf{density weighted payoff} of player~$i$: $w_i(x,t) = u_i(x,t)\cdot q(t)$ for each $x\in X$ and $t\in T$. Let $\cG_i$ be the $\sigma$-algebra generated by the collection of mappings
$$\{w_j(x, \cdot, t_{-i})\colon x\in X, t_{-i} \in T_{-i}, \forall j\neq i\}.$$
Then $\cG_i\subseteq \cT_i$ denotes  player~$i$'s \textbf{inter-player information}. That is, $\cG_i$ is player~$i$'s information flow to all other players, which describes the influence of player~$i$'s private information in other players' payoffs.

\subsection{Existence of pure strategy equilibria}\label{subsec-Bayesian existece}

In this section,  we will prove the  existence of the pure strategy equilibrium in Bayesian games under an appropriate condition called ``coarser inter-player information''. More importantly, we will show that this condition is necessary for the existence result.

\begin{defn}
Player~$i$ is said to have \textbf{coarser inter-player information} if $\cT_i$ has no $\cG_i$-atom under $\lambda_i$.

A Bayesian game is said to have coarser inter-player information if each player has coarser inter-player information.
\end{defn}

\begin{thm}\label{thm-bayesian existence}
Every Bayesian game with coarser inter-player information has a pure strategy equilibrium.
\end{thm}

\begin{rmk}\label{rmk-bayesian}
For Bayesian games with coarser inter-player information, players' payoffs might be interdependent and types could be correlated. In particular, it is inessential whether types are independent or correlated,   since the derivative $q$ can be absorbed into the density weighted payoff.

If $\cG_i$ is the trivial $\sigma$-algebra $\{\emptyset, T_i\}$ for each player $i\in I$, then players have independent priors and private values, and the condition of ``coarser inter-player information'' is automatically satisfied since $(T_i,\cT_i,\lambda_i)$ is atomless.
\end{rmk}

In Theorem~\ref{thm-bayesian existence}, we show that the condition of ``coarser inter-player information'' is sufficient for the existence of pure strategy equilibrium. The next theorem demonstrates that this condition  is also necessary.

Given any $n\ge 2$ and the player space $I = \{1,2,\ldots, n\}$, player~$i$ has private information space $(T_i,\cT_i,\lambda_i)$ and inter-player information $\cG_i$ such that $(T_i,\cG_i,\lambda_i)$ is atomless for each $1\le i\le n$. Let $H_n$ be the collection of all Bayesian games with the player space $I$ and the above private information spaces $\{(T_i,\cT_i/\cG_i,\lambda_i)\}_{i\in I}$.

\begin{thm}\label{thm-bayesian necessity}
Given the player space $I = \{1,\ldots, n\}$ for $n\ge 2$ and the private information space $(T_i,\cT_i/\cG_i,\lambda_i)$ for each $i\in I$, every player~$i$ has coarser inter-player information if either of the following conditions holds:
\begin{enumerate}
  \item every Bayesian game in $H_n$ with type-irrelevant payoffs has a pure strategy equilibrium;\footnote{A Bayesian game is said to have type-irrelevant payoffs if the payoff function of each player does not depend on the type $t\in T$.}
  \item every Bayesian game in $H_n$ with independent types has a pure strategy equilibrium.
\end{enumerate}
\end{thm}

\subsection{Purification}\label{subsec-purification}

In this subsection, we will consider the purification of behavioral strategy profiles in Bayesian games with finite actions.

\begin{defn}\label{defn-undistinguishable}
Let $f=(f_1,f_2,\ldots,f_n)$ and $g=(g_1,g_2,\ldots,g_n)$ be two behavioral strategy profiles.
\begin{enumerate}
  \item The strategy profiles $f$ and $g$ are said to be payoff equivalent if for each player $i \in I$, $U_i(f) = U_i(g)$.

  \item The strategy profiles $f$ and $g$ are said to be strongly  payoff equivalent if
       \begin{enumerate}
         \item they are  payoff equivalent;
         \item for each player $i \in I$ and any given behavioral strategy $h_i$, the two strategy profiles $(h_i, f_{-i})$ and $(h_i, g_{-i})$ are payoff equivalent.
       \end{enumerate}

  \item The strategy profiles $f$ and $g$ are said to be distribution equivalent if for each player $i \in I$, $\int_{T_i} f_i(t_i,\cdot)\rmd \lambda_i(t_i) = \int_{T_i} g_i(t_i,\cdot)\rmd \lambda_i(t_i)$.

  \item Suppose that $f$ is a pure strategy profile. For player $i$, $f_i$ is said to be belief consistent with $g_i$ if $f_i(t_i)\in \supp g_i(t_i)$ for $\lambda_i$-almost all $t_i\in T_i$. Moreover, $f$ is said to be belief consistent with $g$ if they are belief consistent for each player $i\in I$.
\end{enumerate}
\end{defn}

Now we are ready to give the definitions of purification.

\begin{defn}\label{defn-undis}
Suppose that $g$ is a pure strategy profile and $f$ is a behavioral strategy profile. Then $g$ is said to be a strong purification of $f$ if they are strongly payoff equivalent, distribution equivalent, and belief consistent.
\end{defn}

\begin{prop}\label{prop-undis}
In a Bayesian game with coarser inter-player information, every  behavioral strategy profile $f$ possesses a strong purification $g$.
\end{prop}

\section{Appendix}\label{sec-proof}

\subsection{Proofs in Section~\ref{sec-results}}\label{subsec-proof main}

For a sequence of sets $\{A_m\}_{m\in\bN}$ in a topological space $X$, let $\ls(A_m)$ be the set of all $x$ such that for any neighborhood $O_x$ of $x$ there are infinitely many $m$ with $O_x\cap A_m \neq \emptyset$. The following lemma will be needed in the proofs of the main results.

\begin{lem}\label{lem-fatou}
Denote $\{\phi_m\}_{m\in \bN}$ as a sequence of measurable and $p$-integrably bounded mappings from an atomless probability space $(T,\cT,\lambda)$ to $\bR^n$, $1\le p < \infty$.  Let $h_m = E(\phi_m|\cG)$ for each $m\in \bN$, where $\cG$ is a sub-$\sigma$-algebra of $\cT$. Assume that $h_m$ weakly converges to some $h_0 \in L_{p}^\cG(T,\bR^n)$ as $m \to \infty$. If $\cT$ has no $\cG$-atom, then there exists a $\cT$-measurable mapping $\phi_0$ such that
\begin{enumerate}
  \item $\phi_0(t) \in \ls\big( \phi_m(t) \big)$ for $\lambda$-almost all $t\in T$,
  \item $E(\phi_0|\cG) = h_0$.
\end{enumerate}
\end{lem}

\begin{proof}

Since the sequence $\{\phi_m\}_{m\in \bN}$ is $p$-integrably bounded in $L_{p}^\cT(T,\bR^n)$, it has a weakly convergent subsequence by the Riesz/Dunford-Pettis Weak Compactness Theorem in  \citet[p.408/p.412]{RF2010}. Without loss of generality, we assume that  $\phi_m$ weakly converges to some $\phi \in L_{p}^\cT(T,\bR^n)$. Given any $g\in L_{q}^\cG(T,\bR^n)$ such that $\frac{1}{p}+\frac{1}{q} = 1$, we have
\begin{align*}
\int_T h_m g\rmd \lambda
& = \int_T E(\phi_m |\cG) g \rmd \lambda = \int_T E(\phi_m g|\cG) \rmd \lambda \\
& = \int_T \phi_m g \rmd \lambda \to \int_T \phi g \rmd \lambda \\
& = \int_T E(\phi g|\cG) \rmd \lambda = \int_T E(\phi|\cG) g \rmd \lambda.
\end{align*}
Thus, $h_m$ weakly converges to $E(\phi|\cG)$ in $L_{p}^\cG(T,\bR^n)$, which implies that $h_0 = E(\phi|\cG)$. In addition, $\{\phi_k,\phi_{k+1},\ldots\}$ also weakly converges to $\phi$ for each $k\in \bN$. By Theorem~29 of \citet[p.293]{RF2010}, there is a sequence of convex combination of $\{\phi_k,\phi_{k+1},\ldots\}$ that converges to $\phi$ in $L_p$ norm. For each $k\in \bN$, assume that $\varphi_k$ is the convex combination $\{\phi_k,\phi_{k+1},\ldots\}$ such that $\|\varphi_k - \phi\|_{p} \leq \frac{1}{k}$. Thus, there is a subsequence of $\{\varphi_k\}$, say itself, which converges to $\phi$ $\lambda$-almost everywhere.

Fix $t\in T$ such that $\varphi_k(t)$ converges to $\phi(t)$. By Carath\`eodary's convexity theorem (see \citet[Theorem 5.32]{AB2006}),  $\varphi_k (t) = \sum_{j=0}^n \alpha_{jk} \gamma_{jk} (t)$, where
\begin{enumerate}
  \item for each $k\in \bN$, $\alpha_{jk} \geq 0$ for any $j$ and $ \sum_{j=0}^n \alpha_{jk} = 1$;
  \item for each $k\in \bN$, $\gamma_{0k} (t), \ldots, \gamma_{nk} (t) \in \{\phi_k(t),\phi_{k+1}(t),\ldots\}$.
\end{enumerate}
Without loss of generality, assume that for each $0\leq j \leq n$, $\alpha_{jk}\to \alpha_j$ and $\gamma_{jk}(t) \to \gamma_j(t)$. Then $\alpha_{1},\ldots, \alpha_n \geq 0$  and $ \sum_{j=0}^n \alpha_{j} = 1$. Moreover, $\gamma_j(t) \in \ls(\phi_m(t))$. Let $G(t) = \ls(\phi_m(t))$. Then $\phi(t) \in \mbox{co}(G(t))$.

Since $\cT$ has no $\cG$-atom and $G$ is measurable, integrably bounded and closed valued, Theorem~\ref{thm-convex hull} implies that $\cI^{(\cT,\cG)}_G = \cI^{(\cT,\cG)}_{\mbox{co} (G)}$. Thus, there exists a $\cT$-measurable selection $\phi_0$ of $G$ such that $E(\phi_0|\cG) = E(\phi|\cG) = h_0$, which completes the proof.
\end{proof}

\begin{proof}[Proof of Corollary~\ref{coro-convexity}]
First, we assume that $\cT$ has no $\cG$-atom. Pick  two measurable selections  $\phi_1$ and $\phi_2$ of $F$. Let $G(t) = \{\phi_1(t), \phi_2(t)\}$. Then $G$ is a $\cT$-measurable, integrably bounded, closed valued correspondence. By Theorem~\ref{thm-convex hull}, we have $\cI^{(\cT,\cG)}_G = \cI^{(\cT,\cG)}_{\mbox{co} (G)}$, which implies that $\cI^{(\cT,\cG)}_G$ is convex. For any $\alpha\in [0,1]$, there exists a $\cT$-measurable selection $\phi_0$ of $G$ such that $E(\phi_0|\cG) = \alpha E(\phi_1|\cG) + (1-\alpha) E(\phi_2|\cG)$. Since $\phi_0$ is also a selection of $F$, $\cI^{(\cT,\cG)}_F$ is convex.

Conversely, suppose that $\cT$ has a $\cG$-atom $D$ with $\lambda(D) > 0$. Define a correspondence
$$F(t) =
\begin{cases}
\{0,1\} & t\in D;\\
\{0\}   & t\notin D.
\end{cases}$$
It is shown in Proposition~1 of \cite{HS2013} that $\cI^{(\cT,\cG)}_F$ is not convex.\footnote{For simplicity, the target space of the correspondence is $\bR$. One can easily define a new correspondence on $\bR^n$ such that each of other $n-1$ dimensions only contains $0$.}
\end{proof}

Below we prove Theorem~\ref{thm-compact}.

\begin{proof}[Proof of Theorem~\ref{thm-compact}]
Suppose that $\cT$ has no $\cG$-atom and $1 \le p < \infty$. It is sufficient to show that $\cI_{F}^{(\cT,\cG)}$ is weakly sequentially compact in $L_p^\cG(T,\bR^n)$. Fix an arbitrary sequence of $\cT$-measurable selections $\{\phi_m\}_{m\in \bN}$ of $F$. Let $h_m = E(\phi_m|\cG)$ for each $m\in \bN$. We need to show that there is a subsequence of $\{h_m\}_{m\in \bN}$ which weakly converges in $L_p^\cG(T,\bR^n)$ to some point in $\cI_{F}^{(\cT,\cG)}$. Since the sequence $\{\phi_m\}_{m\in \bN}$ is $p$-integrably bounded, it has a weakly convergent subsequence due to the Riesz/Dunford-Pettis Weak Compactness Theorem in  \citet[p.408/p.412]{RF2010}. Without loss of generality, assume that $\phi_m$ weakly converges to some $\phi \in L_{p}^\cT(T,\bR^n)$. As shown in the proof of Lemma~\ref{lem-fatou}, $h_m$ also weakly converges to $E(\phi|\cG)$ in $L_{p}^\cG(T,\bR^n)$. By Lemma~\ref{lem-fatou}, there exists a $\cT$-measurable selection $\phi_0$ of $\ls(\phi_m)$ such that $E(\phi_0|\cG) = E(\phi|\cG)$. Since $F$ is compact valued, $\ls(\phi_m(t)) \subseteq F(t)$ for $\lambda$-almost all $t\in T$. Thus, $\phi_0$ is a selection of $F$, and we are done.

Next, we consider the case $p = \infty$ and $\cG$ is countably generated. Since $F$ is essentially bounded by some positive constant $C$, $\cI_F^{(\cT,\cG)}$ is also norm bounded by $C$. By Alaoglu's Theorem (see Theorem~6.21 of \cite{AB2006}), the closed ball with radius $C$ (the $C$-ball) is weak$^*$ compact in $L_{\infty}^\cG(T,\bR^n)$. We only need to show that $\cI_F^{(\cT,\cG)}$ is weak$^*$ closed in the $C$-ball. Since $\cG$ is countably generated, $L_{1}^\cG(T,\bR^n)$ is separable, which implies that the $C$-ball is metrizable in the weak$^*$ topology (see  Theorem~6.30 of \cite{AB2006}).
Suppose that $\{\phi_m\}$ is a sequence of $\cT$-measurable selections of $F$ and $h_m$ weak$^*$ converges to $h_0 \in L_{\infty}^\cG(T,\bR^n)$  as $m\to \infty$, where $h_m= E(\phi_m|\cG)$ for each $m$.
Then $h_m$ also weakly converges to $h_0$ in $L_{1}^\cG(T,\bR^n)$. Moreover, the condition that $F$ is $\infty$-integrably bounded (i.e., essentially bounded) implies that $F$ is integrably bounded. By Lemma~\ref{lem-fatou}, there exists  a $\cT$-measurable selection $\phi_0$ of $\ls(\phi_m)$ such that $h_0 = E(\phi_0|\cG)$. Since $F$ is compact valued, $\phi_0(t) \in \ls(\phi_m(t)) \subseteq F(t)$ for $\lambda$-almost all $t\in T$. That is, $\phi_0$ is a $\cT$-measurable selection of $F$ and $h_0 \in \cI_F^{(\cT,\cG)}$. Therefore, $\cI_F^{(\cT,\cG)}$ is weak$^*$ closed in the $C$-ball.

\

Conversely, suppose that $\cT$ has a $\cG$-atom $D$ with $\lambda(D) > 0$.
Consider the correspondence $F$ as defined in the proof of Corollary~\ref{coro-convexity}.
Pick an orthonormal subset $\{\varphi_m\}_{m\in \bN}$ of $L_2^{\cT^D}(D,\bR)$ on the atomless probability space $(D,\cT^D,\lambda^D)$ such that $\varphi_m$ takes value in $\{-1,1\}$ and $\int_D \varphi_m\rmd \lambda^D =0$ for each $m\in \bN$. Let
$$\phi_m(t) =
\begin{cases}
\frac{\varphi_m(t)+1}{2} & t\in D;\\
0   & t\notin D.
\end{cases}$$
Then $\phi_m$ is a $\cT$-measurable selection of $F$ for each $m\in \bN$.

Pick a set $E\in \cT^D$. By Bessel's inequality (see \cite[p.316]{RF2010}), $\int_D \mathbf{1}_E \varphi_m\rmd \lambda^D \to 0$ as $m \to \infty$,  where $\mathbf{1}_{E}$ is the indicator function of the set $E$.
Thus, for any $E_1\in \cT$,
\begin{equation} \label{eq-Pi_i}
\int_T \mathbf{1}_{E_1} \phi_m \rmd \lambda  = \frac{1}{2}\int_D \mathbf{1}_{E_1\cap D} \varphi_m \rmd \lambda +\frac{1}{2}\lambda(E_1\cap D) \to \frac{1}{2}\lambda(E_1\cap D).
\end{equation}
Given any nonnegative function $\psi\in L_1^\cT(T,\bR)$, $\psi$ will be the increasing limit of a sequence of simple functions $\{\psi_k\}_{k\in \bN}$ (finite linear combination of measurable indicator functions).
Fix any $\epsilon > 0$. By the dominated convergence theorem, there exists a positive integer $K_0 > 0$ such that for each $k\ge K_0$,
$\int_T |\psi - \psi_k| \rmd \lambda < \epsilon$. Then we have
\begin{align*}
& \quad \left| \int_T \psi \phi_m \rmd \lambda - \frac{1}{2}\int_T  \psi \mathbf{1}_D \rmd \lambda \right| \le \left| \int_T \psi \phi_m \rmd \lambda - \int_T \psi_{K_0} \phi_m \rmd \lambda \right| \\
& + \left| \int_T \psi_{K_0} \phi_m \rmd \lambda - \frac{1}{2}\int_T \psi_{K_0} \mathbf{1}_D \rmd \lambda \right|  + \left| \frac{1}{2}\int_T \psi_{K_0} \mathbf{1}_D \rmd \lambda - \frac{1}{2}\int_T \psi\mathbf{1}_D \rmd \lambda \right| \\
& \le  \int_T | \psi - \psi_{K_0}| \rmd \lambda  + \left| \int_T \psi_{K_0} \phi_m \rmd \lambda - \frac{1}{2}\int_T \psi_{K_0} \mathbf{1}_D \rmd \lambda \right|  + \int_T |\psi_{K_0} - \psi| \rmd \lambda.
\end{align*}
The first and the third terms are less than $\epsilon$. By Equation~(\ref{eq-Pi_i}) and the fact that $\psi_{K_0}$ is a simple function, the second term goes to $0$ as $m \to \infty$. Hence,  $\int_T \psi \phi_m \rmd \lambda \to \frac{1}{2}\int_T  \psi \mathbf{1}_D \rmd \lambda$ as $m \to \infty$. Given any $\psi\in L_1^\cT(T,\bR)$, we can obtain $\int_T \psi \phi_m \rmd \lambda \to \frac{1}{2}\int_T  \psi \mathbf{1}_D \rmd \lambda$ as $m \to \infty$ by writing $\psi$ as the sum of its positive and negative parts.
Therefore, $\phi_m$ weak$^*$ converges to $\phi = \frac{1}{2}\mathbf{1}_D$ in $L_\infty^\cT(T,\bR)$. Thus, $E(\phi_m|\cG) \in \cI^{(\cT,\cG)}_F$ weak$^*$ converges to $\frac{1}{2} E(\mathbf{1}_D|\cG)$ in $L_\infty^\cG(T,\bR)$ as shown in the proof of Lemma~\ref{lem-fatou}. It is shown in \cite{HS2013} that $\frac{1}{2} E(\mathbf{1}_D|\cG) \notin \cI^{(\cT,\cG)}_F$, which implies that $\cI^{(\cT,\cG)}_F$ is not weak$^*$ compact in $L_\infty^\cG(T,\bR)$.

For $1\le p < \infty$, just note that $F$ is also $p$-integrably bounded, and $\phi_m$ weakly  converges to $\phi = \frac{1}{2}\mathbf{1}_D$ in $L_{p}^\cT(T,\bR)$.
\end{proof}

\begin{proof}[Proof of Theorem~\ref{thm-upper}]
Suppose that $\cT$ has no $\cG$-atom and  $1\le p < \infty$.
By Theorem~\ref{thm-compact}, we know that $\cI_G^{(\cT,\cG)}$ is weakly compact, and hence weakly sequentially compact.
Pick $\{y_m\}_{m=0}^\infty \subseteq Y$ and $\{\phi_m\}_{m\in\bN}$ such that $\phi_m$ is a $\cT$-measurable selection of $F_{y_m}$.
Let $h_m = E(\phi_m|\cG)$ for each $m\in \bN$. Suppose that $h_m$ weakly converges to some $h_0\in L_{p}^\cG(T,\bR^n)$ and $y_m$ converges to some $y_0\in Y$.
By Lemma~\ref{lem-fatou}, there exists  a $\cT$-measurable selection $\phi_0$ of $\ls(\phi_m)$ such that $h_0 = E(\phi_0|\cG)$.

Since $F_t(\cdot)$ is upper hemicontinuous for $\lambda$-almost all $t\in T$, $\phi_0(t) \in \ls(\phi_m(t)) \subseteq \ls( F_{y_m}(t)) \subseteq F_{y_0}(t)$ for $\lambda$-almost all $t\in T$.
That is, $\phi_0$ is a $\cT$-measurable selection of $F_{y_0}$ and $h_0 \in H(y_0)$. Therefore, $H$ is weakly upper hemicontinuous.
The case that $p = \infty$ and  $\cG$ is countably generated follows from a similar argument by noting that any closed ball in $L_\infty^\cG(T,\bR^n)$ is metrizable.

\

Conversely, suppose that $\cT$ has a $\cG$-atom $D$ with $\lambda(D) > 0$. Let $G$ be the correspondence as in Corollary~\ref{coro-convexity}
$$G(t) =
\begin{cases}
\{0,1\} & t\in D;\\
\{0\}   & t\notin D.
\end{cases}$$
Let $Y = \{\frac{1}{m}\}_{m \ge 1} \cup \{0\}$ endowed with the usual metric, $F(t,0) = G(t)$ and $F(t, \frac{1}{m}) = \{\phi_m(t)\}$ for all $t\in T$ and $m\ge 1$, where $\phi_m$ is the same as in the converse part of the proof of Theorem~\ref{thm-compact}.
Then $G$ is compact valued and bounded, and $F(t,\cdot)$ is upper hemicontinuous for all $t\in T$.

Consider the correspondence $G$. For $1\le p < \infty$, since $\cI_G^{(\cT,\cG)}$ is $p$-integrably bounded, it is relatively weakly sequentially compact in $L_{p}^\cG(T,\bR)$ due to the Riesz/Dunford-Pettis Weak Compactness Theorem in  \citet[p.408/p.412]{RF2010}, and hence relatively weakly compact. For $p = \infty$, $\cI_G^{(\cT,\cG)}$ is  relatively weak$^*$ compact in $L_{\infty}^\cG(T,\bR^n)$ due to Alaoglu's Theorem. Thus, $H(y)$ is a subset of a fixed weakly (resp. weak$^*$) compact set  for all $y \in Y$ when $1\le p < \infty$ (resp. $p=\infty$).

For the sequence $\{\frac{1}{m}\}$, $\frac{1}{m} \to 0$ and $\phi_m$ is a selection of $F_{\frac{1}{m}}$. As shown in the proof above, $E(\phi_m|\cG)$ weakly (resp. weak$^*$) converges to $\frac{1}{2} E(\mathbf{1}_D|\cG)$ in $L_{p}^\cG(T,\bR)$ for $1 \le p < \infty$ (resp. $p = \infty$), but there is no $\cT$-measurable selection $\phi_0$ of $G$ such that $E(\phi_0|\cG) = \frac{1}{2} E(\mathbf{1}_D|\cG)$. Therefore, $\frac{1}{2} E(\mathbf{1}_D|\cG) \notin H(0)= \cI_{F_0}^{(\cT,\cG)}$, which implies that $H(y)$ is neither weakly upper hemicontinuous in $L_{p}^\cG(T,\bR)$  for $1 \le p < \infty$ nor weak$^*$ upper hemicontinuous in $L_{\infty}^\cG(T,\bR)$.
\end{proof}

\subsection{Proofs in Section~\ref{sec-Bayesian game}}

\subsubsection{Proofs in Section~\ref{subsec-Bayesian existece}}

\begin{proof}[Proof of Theorem~\ref{thm-bayesian existence}]
Suppose that $m_i$ is the cardinality of $X_i$ for each $i\in I$. Let $\triangle_i = \{(a^i_1,\ldots, a^i_{m_i}) \colon \sum_{1\le k \le m_i} a^i_k = 1, a^i_k \ge 0 \mbox{ for } 1\le k \le m_i\}$. The mixture of actions of player~$i$ in $\cM(X_i)$ can be regarded as elements in the simplex $\triangle_i$, and the pure actions in $X_i$ correspond to the extreme points of $\triangle_i$. Denote $L_i^{\cG_i}$ as the set of all $\cG_i$-measurable functions from $T_i$ to $\cM(X_i)$.  Without loss of generality, it can be viewed as $L_{\infty}^{\cG_i}(T_i,\triangle_i)$, and embedded in $L_{\infty}^{\cG_i}(T_i,\bR^{m_i})$ endowed with the weak$^*$ topology.

By the Riesz representation theorem (see Theorem~13.28 of \cite{AB2006}), $L^{\cG_i}_{\infty}(T_i, \bR^{m_i})$ can be viewed as the dual space of $L_{1}^{\cG_i}(T_i, \bR^{m_i})$. Then $L^{\cG_i}_{\infty}(T_i, \bR^{m_i})$ is a locally convex, Hausdorff topological vector space under the weak$^*$ topology.   By Alaoglu's Theorem (see Theorem~6.21 of \cite{AB2006}), the closed ball with radius $C \ge 1$ (the $C$-ball) is weak$^*$ compact in $L^{\cG_i}_{\infty}(T_i, \bR^{m_i})$. Since $\cG_i$ is countable generated, $L^{\cG_i}_1(T_i, \bR^{m_i})$ is separable, which implies that the $C$-ball is metrizable in the weak$^*$ topology (see Theorem~6.30 of \cite{AB2006}).

It is obvious that $L_i^{\cG_i}$ is included in the $C$-ball, we need to show that $L_i^{\cG_i}$ is weak$^*$ closed. That is, for any sequence $\{g_k\}_{k\in \bN} \in L_i^{\cG_i}$ such that $g_k$ weak$^*$ converges to some $g_0 \in L^{\cG_i}_{\infty}(T_i, \bR^{m_i})$, we need to show $g_0 \in L_i^{\cG_i}$. Since $g_k$ weak$^*$ converges to $g_0 \in L^{\cG_i}_{\infty}(T_i, \bR^{m_i})$, it also weakly converges to $g_0$  in $L^{\cG_i}_{1}(T_i, \bR^{m_i})$. Following an analogous argument in the proof of Lemma~\ref{lem-fatou}, one can show that $g_0(t_i) \in \mbox{co}\left(\ls(g_k(t_i))\right)$ for $\lambda_i$-almost all $t_i \in T_i$. Since $\triangle_i$ is closed and convex, $\mbox{co}\left(\ls(g_k(t_i))\right) \subseteq \triangle_i$ for $\lambda_i$-almost all $t_i \in T_i$. Thus, $g_0(t_i) \in \triangle_i$ for $\lambda_i$-almost all $t_i \in T_i$, and $g_0 \in L_i^{\cG_i}$. Therefore, $L_i^{\cG_i}$ is nonempty, convex, and compact under the weak$^*$ topology. Let $L^{\cG} = \times_{i\in I} L_i^{\cG_i}$ endowed with the product topology.

Given a pure strategy profile $h$, let $\overline{h_i} = E^{\lambda_i}(h_i|\cG_i) \in L_i^{\cG_i}$, and $h_i^k$ denote the $k$-th dimension of $h_i$ for each player~$i \in I$ and $1\le k \le m_i$. For any distinct $i,j \in I$, $x_{-j} \in X_{-j}$, $t_{-j} \in T_{-j}$, and $D\in \cG_j$,
\begin{align*}
& \quad \int_{T_{j}}  \mathbf{1}_D(t_j) u_i(x_{-j},h_{j}(t_{j}), t_{-j}, t_j) q(t_j, t_{-j})  \lambda_{j} (\rmd t_{j}) \\
& = \int_{T_{j}}  \mathbf{1}_D(t_j) w_i(x_{-j},h_{j}(t_{j}), t_{-j}, t_j)  \lambda_{j} (\rmd t_{j}) \\
& = \int_{T_{j}}  \mathbf{1}_D(t_j) \sum_{k=1}^{m_j} \left(w_i(x_{-j},a^j_{k}, t_{-j}, t_j) \cdot h_{j}^k(t_{j})\right)  \lambda_{j} (\rmd t_{j}) \\
& = \int_{T_{j}}  E^{\lambda_{j}}\left(\mathbf{1}_D(t_j) \sum_{k=1}^{m_j} \left(w_i(x_{-j},a^j_{k}, t_{-j}, t_j) \cdot h_{j}^k(t_{j})\right) |\cG_{j}\right)  \lambda_{j} (\rmd t_{j}) \\
& = \int_{T_{j}}  \mathbf{1}_D(t_j) \sum_{k=1}^{m_j}E^{\lambda_{j}}\left(w_i(x_{-j},a^j_{k}, t_{-j}, t_j) \cdot h_{j}^k(t_{j})|\cG_{j}\right) \lambda_{j} (\rmd t_{j}) \\
& = \int_{T_{j}}  \mathbf{1}_D(t_j) \sum_{k=1}^{m_j} w_i(x_{-j},a^j_{k}, t_{-j}, t_j) \cdot E^{\lambda_{j}}\left( h_{j}^k(t_{j})|\cG_{j}\right) \lambda_{j} (\rmd t_{j}) \\
& = \int_{T_{j}}  \mathbf{1}_D(t_j) \sum_{k=1}^{m_j} w_i(x_{-j},a^j_{k}, t_{-j}, t_j) \cdot \overline{h_j^k}(t_{j}) \lambda_{j} (\rmd t_{j}) \\
& = \int_{T_{j}}  \mathbf{1}_D(t_j) \int_{X_{j}} w_i(x_{-j},x_{j}, t_{-j}, t_j) \overline{h}_{j} (t_{j}, \rmd x_{j})  \lambda_{j} (\rmd t_{j}) .
\end{align*}
The first three equalities are obvious. The  fourth and fifth equalities hold since $\mathbf{1}_{D}$ and $w_i(x, t_{-j}, \cdot)$ are $\cG_j$-measurable for any $x\in X$ and $t_{-j} \in T_{-j}$. The sixth equality is due to the definition of $\overline{h_j}$, and the last one is just rewriting the summation as integration. Thus, for $\lambda_j$-almost all $t_j \in T_j$,
\begin{equation}\label{equa-3}
E^{\lambda_{j}}\left(w_i(x_{-j},h_{j}(t_{j}), t_{-j}, t_j) |\cG_{j}\right) = \int_{X_{j}} w_i(x_{-j},x_{j}, t_{-j}, t_j) \overline{h}_{j} (t_{j}, \rmd x_{j}).
\end{equation}

Fix player~$1$. For any $t_1\in T_1$ and $x_1\in X_1$, we have
\begin{align*}
& \quad \int_{T_{-1}}  u_1(x_1,h_{-1}(t_{-1}), t_1,t_{-1})  q(t_1,t_{-1})  \lambda_{-1} (\rmd t_{-1}) \\
& = \int_{T_{-1}}  w_1(x_1,h_{-1}(t_{-1}), t_1,t_{-1})  \lambda_{-1} (\rmd t_{-1}) \\
& = \int_{T_{-(1,2)}}  \int_{T_2} w_1(x_1, h_{2}(t_{2}), h_{-(1,2)}(t_{-(1,2)}), t_{-2},t_{2})  \lambda_{2} (\rmd t_{2}) \lambda_{-(1,2)} (\rmd t_{-(1,2)}) \\
& = \int_{T_{-(1,2)}}  \int_{T_2} E^{\lambda_{2}}\left(w_1(x_1, h_{2}(t_{2}), h_{-(1,2)}(t_{-(1,2)}), t_{-2},t_{2}) |\cG_{2}\right) \lambda_{2} (\rmd t_{2}) \lambda_{-(1,2)} (\rmd t_{-(1,2)}) \\
& = \int_{T_{-(1,2)}}  \int_{T_2} \int_{X_{2}} w_1(x_1, x_{2}, h_{-(1,2)}(t_{-(1,2)}),  t_{-2},t_{2}) \overline{h}_{2} (t_{2}, \rmd x_{2}) \lambda_{2} (\rmd t_{2}) \lambda_{-(1,2)} (\rmd t_{-(1,2)}) \\
& = \cdots \\
& = \int_{T_{-1}}  \int_{X_{-1}} w_1(x_1, x_{-1}, t_1,t_{-1}) \overline{h}_{-1} (t_{-1}, \rmd x_{-1}) \lambda_{-1} (\rmd t_{-1}),
\end{align*}
where the subscript $-(1,2)$ denotes the set of all players except players~$1$ and $2$. The first equality is due to the definition of density weighted payoff. The second equality is due to the Fubini property. The third equality holds by taking the conditional expectation. The fourth equality is implied by Equation~(\ref{equa-3}). Then the previous four equalities are repeated for $n-2$ times (from $T_3$  to $T_n$). This procedure is omitted in the fifth equality, and finally leads to the last equality. One can repeat the argument and show that for any $i\in I$, $x_i \in X_i$ and $t_i \in T_i$
\begin{equation}\label{equa-4}
\begin{split}
& \quad  \int_{T_{-i}}  u_i(x_i,h_{-i}(t_{-i}), t_i,t_{-i})  q(t_i,t_{-i})  \lambda_{-i} (\rmd t_{-i}) \\
& = \int_{T_{-i}}  \int_{X_{-i}} w_i(x_i, x_{-i}, t_i,t_{-i}) \overline{h}_{-i} (t_{-i}, \rmd x_{-i}) \lambda_{-i} (\rmd t_{-i}).
\end{split}
\end{equation}

For each $i\in I$, let $F_i$ be a mapping from $T_i\times X_i\times L^{\cG}$ to $\bR$ defined as follows:
$$F_i(t_i,x_i,g_1,\ldots, g_n) = \int_{T_{-i}}  \int_{X_{-i}} w_i(x_i, x_{-i}, t_i,t_{-i})  g_{-i} (t_{-i}, \rmd x_{-i}) \lambda_{-i} (\rmd t_{-i}).
$$
It is clear that $F_i$ is $\cT_i$-measurable on $T_i$ and continuous on $L^{\cG}$, where $L^{\cG}$ is endowed with the weak$^*$ topology. For each $i\in I$, the best response correspondence $G_i$ from $T_i\times L^{\cG}$ to $X_i$ is given by
$$G_i(t_i,g_1,\ldots, g_n) = \mbox{argmax}_{x_i\in X_i} F_i(t_i,x_i,g_1,\ldots, g_n).$$
For each $t_i$, Berge's maximal theorem implies that $G_i$ is nonempty, compact-valued, and upper-hemicontinuous on $L^{\cG}$. For any $x_i$ and $(g_1,\ldots, g_n)$, $F_i$ is  $\cT_i$-measurable. Then  $G_i(\cdot,g_1,\ldots,g_n)$ admits a $\cT_i$-measurable selection. Thus, $E^{\lambda_i}\left(G_i(\cdot,g_1,\ldots,g_n)|\cG_i\right)$ is nonempty. Since $\cT_i$ has no $\cG_i$-atom, by Corollary~\ref{coro-convexity} and Theorems~\ref{thm-compact} and \ref{thm-upper},  it is convex, weak$^*$ compact-valued, and weak$^*$ upper-hemicontinuous on $L^{\cG}$.

Consider a correspondence from $L^{\cG}$ to itself:
$$\psi(g_1,\ldots,g_n)=\times_{i=1}^{n}E^{\lambda_i}\left(G_i(\cdot,g_1,\ldots,g_n)|\cG_i\right).$$
It is clear that $\psi$ is nonempty, convex, weak$^*$ compact-valued, and weak$^*$ upper-hemicontinuous on $L^{\cG}$. By Fan-Glicksberg's fixed-point theorem, there exists a fixed point $(g^*_1,\ldots,g^*_n)$ of $\psi$. Thus for each $i$, there exists some $\cT_i$-measurable selection $f_i^*$  of $G_i(\cdot,g_1^*,\ldots,g_n^*)$ such that $g_i^* = E^{\lambda_i}\left(f_i^*|\cG_i\right)$.

With the strategy profile $(f_1^*,\ldots,f_n^*)$, the payoff of player $i$ is
\begin{align*}
U_i(f^*)
& = \int_T w_i(f_i^*(t_i),f_{-i}^*(t_{-i}),t_i,t_{-i}) \lambda(\rmd t) \\
& = \int_{T_i} \int_{T_{-i}} w_i(f_i^*(t_i),f_{-i}^*(t_{-i}),t_i, t_{-i})  \lambda_{-i}(\rmd t_{-i})  \lambda_{i}(\rmd t_i) \\
& = \int_{T_i} \int_{T_{-i}}  \int_{X_{-i}} w_i(f_i^*(t_i), x_{-i}, t_i,t_{-i}) g_i^* (t_{-i}, \rmd x_{-i}) \lambda_{-i} (\rmd t_{-i})  \lambda_{i}(\rmd t_i).
\end{align*}

The first equality holds due to the definition of $U_i$. The second equality holds based on  the Fubini property, and the third equality relies on Equation~(\ref{equa-4}). By the choice of $(g_1^*,\ldots,g_n^*)$, we have that $(f_1^*,\ldots,f_n^*)$ is a pure strategy equilibrium.
\end{proof}

To prove Theorem~\ref{thm-bayesian necessity}, we first consider an auxiliary game.

\begin{exam}\label{exma-1}
Consider an $m\times m$ zero-sum ``matching pennies'' game $\Gamma$ with asymmetric information. There are two players, and the action space for both players is $A_1 = A_2 =\{a_1,a_2,\ldots a_m\}$, $m\ge 2$. The payoff matrix for player 1 is given below.
\begin{figure}[ht]\hspace*{\fill}%
    \begin{game}{5}{5}[Player $1$][Player $2$][]
        		        & $a_1$	    & $a_2$         & $a_3$	        & $\cdots$		& $a_m$		\\
        $a_1$	        & $1$	    & $-1$	        & $0$	        & $\cdots$      & $0$	\\
        $a_2$	        & $0$	    & $1$	        & $-1$	        & $\cdots$      & $0$	\\
        $a_3$	        & $0$	    & $0$	        & $1$	        & $\cdots$      & $0$	\\
        $\vdots$	& $\vdots$	& $\vdots$	    & $\vdots$	    & $\vdots$	    & $\vdots$	\\
        $a_m$	        & $-1$	    & $0$           & $\cdots$	    & $0$	        & $1$	\\
     \end{game}\hspace*{\fill}%
\end{figure}

Player $i$ has a private information space $L_i = [0,1]$ and  $(l_1,l_2)$ follows the uniform distribution $\tau$ on the triangle of the unit square $0 \le l_1 \le l_2 \le 1$. Then it is obvious that the Radon-Nikodym derivative  of $\tau$ with respect to the Lebesgue measure on the unit square is
$$ \rho(l_1, l_2) =
\begin{cases}
2, & 0 \le l_1 \le l_2 \le 1; \\
0, & \mbox{otherwise.}
\end{cases}
$$

Let $\tau_i$ be the marginal distribution of $\tau$ on $L_i$ for $i = 1,2$. Then the Lebesgue measure $\eta$ is absolutely continuous with respect to $\tau_1$  on $[0,1]$ with the Radon-Nikodym derivative $\beta_1(l_1) = \frac{1}{2(1-l_1)}$ if $0 < l_1 < 1$ and $0$ otherwise, and $\eta$ is absolutely continuous with respect to  $\tau_2$ with the Radon-Nikodym derivative $\beta_2(l_2) = \frac{1}{2l_2}$ if $0 < l_2 < 1$ and $0$ otherwise. The two probability measures $\tau_1$ and $\tau_2$ are both atomless. Let $\rho'$ be the corresponding Radon-Nikodym derivative  of $\tau$ with respect to $\tau_1\otimes \tau_2$:
$$ \rho'(l_1, l_2) = \rho(l_1, l_2)\cdot \beta_1(l_1) \cdot \beta_2(l_2) =
\begin{cases}
\frac{1}{2(1-l_1)l_2}, & 0 < l_1 \le l_2 < 1; \\
0, & \mbox{otherwise.}
\end{cases}
$$
\end{exam}

As is well known, there exists a measure preserving mapping $h_i$ from $(T_i, \cG_i, \lambda_i)$ to  $([0,1],\cB,\tau_i)$ such that for any $E\in\cG_i$, there exists a set $E'\in \cB$ such that $\lambda_i(E\triangle h_i^{-1}(E'))=0$. For $i = 1,2$, let $\pi_i$ be a probability measure on $(T_i, \cT_i)$ which is absolutely continuous with respect to $\lambda_i$ with the Radon-Nikodym derivative $\beta_i(h_i(t_i))$. Since $\beta_i(h_i(t_i))$ is positive for $\lambda_i$-almost all $t_i$, $\lambda_i$ is also absolutely continuous with respect to $\pi_i$.

\begin{proof}[Proof of Theorem~\ref{thm-bayesian necessity}]
\

(1) First we consider the following $2$-player game $\Gamma'$, and then extend it to an $n$-player game.  Player~$1$ and $2$'s action spaces and payoffs are the same as in the game $\Gamma$. The private information space for player $i$ is $(T_i, \cT_i, \lambda_i)$, $q(t_1, t_2) = \rho'(h_1(t_1), h_2(t_2))$, and the common prior $\lambda$ has the Radon-Nikodym derivative $q$ with respect to $\lambda_1 \otimes \lambda_2$. It can be easily checked that $\lambda_i$ is the marginal probability measure of $\lambda$ for $i = 1,2$.

Suppose that $\Gamma_1$ has a pure strategy equilibrium $(f_1, f_2)$. Let $E^1_j = \{t_1\in T_1\colon f_1(t_1) = a_j\}$ and $E^2_j = \{t_2\in T_2\colon f_2(t_2) = a_j\}$ for $1\le j \le m$. Then we shall show that for $\lambda_2$-almost all $t_2\in T_2$
\begin{equation}\label{equa-1}
\int_{E^1_{1}} q(t_1, t_2) \lambda_1(\rmd t_1) = \ldots = \int_{E^1_{m}} q(t_1, t_2) \lambda_1(\rmd t_1),
\end{equation}
and for $\lambda_1$-almost all $t_1\in T_1$,
\begin{equation}\label{equa-2}
\int_{E^2_{1}} q(t_1, t_2) \lambda_2(\rmd t_2) = \ldots = \int_{E^2_{m}} q(t_1, t_2) \lambda_2(\rmd t_2).
\end{equation}

Suppose that $\alpha$ and  $-\alpha$ are the equilibrium payoffs of player~1 and player~2, respectively. Denote $a_{m+1} = a_1$,  $E^1_{m+1} = E^1_{1}$ and $E^2_{m+1} = E^2_{1}$.  For $j = 1, \ldots, m$, let
$$C^1_j = \left\{t_1\in T_1 \colon  \int_{E^2_{j}} q(t_1, t_2) \lambda_2(\rmd t_2) > \int_{E^2_{j+1}} q(t_1, t_2) \lambda_2(\rmd t_2)  \right\}
$$
and
$$C^2_j = \left\{t_2\in T_2 \colon  \int_{E^1_{j}} q(t_1, t_2) \lambda_1(\rmd t_1) > \int_{E^1_{j+1}} q(t_1, t_2) \lambda_1(\rmd t_1)  \right\}.
$$

Now we define a new strategy for players~1 and 2 as follows:
$$ f_1'(t_1) =
\begin{cases}
a_{j} & t_1\in C^1_j \setminus \left( \cup_{1\le k < j} C^1_k \right), \\
a_1     & \mbox{otherwise};
\end{cases}
$$
and
$$ f_2'(t_2) =
\begin{cases}
a_{j+1} & t_2\in C^2_j \setminus \left( \cup_{1\le k < j} C^2_k \right), \\
a_1     & \mbox{otherwise}.
\end{cases}
$$
We claim that player~2 can choose the strategy $f_2'$ and get a nonnegative payoff.

If player~$2$ takes action $a_{j+1}$ at state $t_2$, then his interim expected payoff is
\begin{align*}
\quad \int_{T_1} u_2(f_1(t_1), a_{j+1}) q(t_1, t_2) \lambda_1(t_1)
& = \sum_{k=1}^m \int_{E^1_k} u_2(a_k, a_{j+1}) q(t_1, t_2) \lambda_1(t_1) \\
& = \int_{E^1_{j}} q(t_1, t_2) \lambda_1(t_1) -\int_{E^1_{j+1}}  q(t_1, t_2) \lambda_1(t_1).
\end{align*}
\begin{enumerate}
  \item if $t_2 \in C^2_j$, then choosing the action $a_{j+1}$ gives player~2 a  strictly positive payoff;
  \item if $t_2 \in T_2\setminus \left(\cup_{1\le j \le m} C^2_j \right)$, then player~2 is indifferent between any action and gets a payoff $0$.
\end{enumerate}
Thus, player~2 can choose the strategy $f_2'$ and guarantee himself a nonnegative payoff, which implies that $\alpha \le 0$. Similarly, one can analyze the payoff of player~1 and show that $\alpha \ge 0$. As a result, $\alpha = 0$, which implies that $\lambda_1\left(\cup_{1\le j \le m} C^1_j \right) = 0$ and $\lambda_2\left(\cup_{1\le j \le m} C^2_j \right) = 0$. As a result, Equation~(\ref{equa-1}) holds for $\lambda_2$-almost all $t_2\in T_2$, and  Equation~(\ref{equa-2}) holds for $\lambda_1$-almost all $t_1\in T_1$.

For $\lambda_2$-almost all $t_2\in T_2$, $\int_{T_1} q(t_1, t_2) \lambda_1(\rmd t_1) = 1$, we have
$$\int_{E^1_{j}} q(t_1, t_2) \lambda_1(\rmd t_1) = \frac{1}{m}$$
for  $1\le j\le m$. For each $E^1_{j}$, there exists a set $D_j \subseteq L_1$ such that $\lambda_1(E_j^1 \triangle h_1^{-1}(D_j)) = 0$, implying that  $$\int_{D_{j}} \rho'(l_1, h_2(t_2)) \tau_1(\rmd l_1) = \int_{D_{j}} \rho(l_1, h_2(t_2)) \beta_2(h_2(t_2))\eta(\rmd l_1) = \frac{1}{m},$$
and hence $\int_{D_{j}} \rho(l_1, h_2(t_2))\eta(\rmd l_1) = \frac{ 2(h_2(t_2))}{m}$ for $\lambda_2$-almost all $t_2\in T_2$.  That is, for $\lambda_2$-almost all $t_2\in T_2$,  $\eta(D_j \cap [0, h_2(t_2)]) = \frac{h_2(t_2)}{m}$.

As a result,
\begin{align*}
\pi_1(E^1_j \cap h_1^{-1}\left([0, h_2(t_2)]\right))
& = \int_{E^1_j \cap h_1^{-1}\left([0, h_2(t_2)]\right)} \beta_1(h_1(t_1)) \lambda_1(\rmd t_1) \\
& = \int_{D_j \cap [0, h_2(t_2)]} \beta_1(l_1) \tau_1(\rmd l_1) \\
& = \eta(D_j \cap [0, h_2(t_2)]) \\
& = \frac{h_2(t_2)}{m}.
\end{align*}
Thus, $\pi_1(E^1_j) = \frac{1}{m}$. In addition, $\pi_1(h_1^{-1}\left([0, h_2(t_2)]\right)) = \eta \left([0, h_2(t_2)]\right) = h_2(t_2)$ for $\lambda_2$-almost all $t_2 \in T_2$. Therefore, $\pi_1(E^1_j \cap h_1^{-1}\left([0, h_2(t_2)]\right)) = \pi_1(E^1_j) \cdot \pi_1(h_1^{-1}\left([0, h_2(t_2)]\right))$ for $\lambda_2$-almost all $t_2\in T_2$. Since $\{[0, h_2(t_2)]\}_{t_2 \in T_2}$ generates the Borel $\sigma$-algebra on $[0,1]$ modulo null sets, $\{h_1^{-1}\left([0, h_2(t_2)]\right)\}_{t_2 \in T_2}$  generates $\cG_1$ on $T_1$ modulo null sets, which implies that $E_j^1$ is independent of $\cG_1$ under $\pi_1$. As $m$ is arbitrary, we have proved that for any natural number $m\ge 2$, there exist $m$ disjoint subsets $\{E_j^1\}_{1\le j \le m}$ which are of measure $\frac{1}{m}$ and independent of $\cG_1$ under $\pi_1$. Thus, $\cT_1$ has no $\cG_1$-atom under $\pi_1$. Since $\pi_1$ and $\lambda_1$ are absolutely continuous with respect to each other, $\cT_1$ has no $\cG_1$-atom under $\lambda_1$.  Similarly, one can show that $\cT_2$ has no $\cG_2$-atom under $\lambda_2$.

We extend the game $\Gamma'$ to an $n$-player game $\Gamma_2$. Players~$1$ and $2$ in $\Gamma_2$ share the same payoffs, action sets and private information spaces with those in the game $\Gamma'$. Other players in $\Gamma_2$ are dummy in the sense that player~$k$ has private information space $(T_k,\cT_k,\lambda_k)$, and only one action set $X_k = \{a\}$ for $3 \le k \le n$. The common prior $\lambda$ is absolutely continuous with respect to $\bigotimes_{1\le i \le n}\lambda_i$  with the Radon-Nikodym derivative $q(t_1,t_2)$. Hence, the  payoffs of all players in the Bayesian game $\Gamma_2$  are type-irrelevant. If $\Gamma_2$ has a pure strategy equilibrium, then the analysis above shows that players~$1$ and $2$ have coarser inter-player information.

For any $3\le j \le n$, one can  construct a new $n$-player game $\Gamma_j$ in which players~$1$ and $j$ are active while all other players are dummy. The payoff functions,  action sets and private information spaces of players~$1$ and $j$ are defined similarly as those of players~$1$ and $2$ in the game $\Gamma_2$. Adopting the above argument, it can be shown that players~$1$ and $j$ have coarser inter-player information. Therefore,  all players have coarser inter-player information.

\

(2) Now we construct a new game $\Gamma''$ based on the game $\Gamma'$ above. Suppose that players~$1$ and $2$'s action spaces and private information spaces are the same, while the payoff of player $i$ in the game $\Gamma''$ is given by $v_i(x,t) = u_i(x)\cdot q(t)$ for each $x\in X$ and $t\in T$, where $u_i$ is the payoff function of player~$i$ and $q$ is the Radon-Nikodym derivative in the game $\Gamma'$. Players have independent types and the common prior $\lambda = \lambda_1 \otimes \lambda_2$. It is obvious that the game $\Gamma''$ is essentially the same compared with the game $\Gamma'$ if one considers the density weighted payoff. Extending the game $\Gamma''$ to an $n$-player game with $2$ active players and $n-2$ dummy players,  then one can follow the proof in (1) and show that the two active players have coarser inter-player information if the game $\Gamma''$ has a pure strategy equilibrium. Since those two active players are arbitrarily chosen, all players have coarser inter-player information.
\end{proof}

\subsubsection{Proofs in Section~\ref{subsec-purification}}

\begin{proof}[Proof of Proposition~\ref{prop-undis}]
Given any behavioral strategy profile $f$, $x_i \in X_i$ and $t_i \in T_i$, let
$$V_i^{f}(x_i,t_i) = \int_{T_{-i}} \int_{X_{-i}} u_i(x_i,x_{-i}, t_i,t_{-i})  q(t_i,t_{-i})  \prod_{j \neq i} f_j(t_j,\rmd x_j) \lambda_{-i} (\rmd t_{-i}).
$$
For any  $\mu_i\in \cM(X_i)$, define
$$W_i^{(\mu_i,f)}(t_i) = \int_{X_i} V_i^f(x_i,t_i) \mu_i(\rmd x_i).
$$
Let $c_i(t_i,x_i)=\mathbf{1}_{\supp f_i(t_i)}(x_i)$ for each $t_i\in T_i$, $x_i\in X_i$ and $i\in I$. Denote $c_i^{\mu_i}(t_i) = \int_{X_i} c_i(t_i,x_i) \mu_i(\rmd x_i)$ for any $\mu_i \in \cM(X_i)$. Then given any behavioral strategy $h_i$, we slightly abuse the notation by letting $W_i^{(h_i,f)}(t_i) = \int_{X_i} V_i^f(x_i,t_i) h_i(t_i,\rmd x_i)$ and $c_i^{h_i}(t_i) = \int_{X_i} c_i(t_i,x_i) h_i(t_i,\rmd x_i)$.

Define a correspondence
$$H_i^{f}(t_i) = \{\left(x_i, V_i^{f}(x_i,t_i), c_i(t_i,x_i) \right) \colon x_i\in X_i \}.$$
We have
$$co(H_i^{f})(t_i) = \{\left(\mu_i, W_i^{(\mu_i,f)}(t_i), c_i^{\mu_i}(t_i) \right) \colon \mu_i\in \cM(X_i) \}.$$
By Theorem~\ref{thm-convex hull}, we have $E^{\lambda_i}(H_i^f|\cG_i) = E^{\lambda_i}(co(H_i^{f}) |\cG_i)$.

For each $i\in I$,  $(f_i, W_i^{(f_i,f)},  c_i^{f_i})$ is a measurable selection of $co(H_i^{f})$. Thus, there is a $\cT_i$-measurable mapping $g_i$ from $T_i$ to $X_i$ such that $E^{\lambda_i}(g_i|\cG_i) = E^{\lambda_i}(f_i|\cG_i)$, $E^{\lambda_i}( W_i^{(g_i,f)}|\cG_i) = E^{\lambda_i}(W_i^{(f_i,f)}|\cG_i)$ and $E^{\lambda_i}(c_i^{g_i}|\cG_i) = E^{\lambda_i}(c_i^{f_i}|\cG_i)$. Then $E^{\lambda_i}(g_i|\cG_i) = E^{\lambda_i}(f_i|\cG_i)$ for each $i$ implies that $f$ and $g$ are distribution equivalent.

Given any $t_i\in T_i$ and $x_i\in X_i$,
\begin{align*}
V_i^{g}(x_i,t_i)
& = \int_{T_{-i}} \int_{X_{-i}} u_i(x_i,x_{-i}, t_i,t_{-i})  q(t_i,t_{-i})  \prod_{j \neq i} g_j(t_j,\rmd x_j) \lambda_{-i} (\rmd t_{-i})\\
& = \int_{T_{-i}} \int_{X_{-i}} w_i(x_i,x_{-i}, t_i,t_{-i})  \prod_{j \neq i} g_j(t_j,\rmd x_j) \lambda_{-i} (\rmd t_{-i})\\
& = \int_{T_{-i}} \int_{X_{-i}} w_i(x_i,x_{-i}, t_i,t_{-i})  \prod_{j \neq i} E^{\lambda_j}\left(g_j|\cG_j \right)(t_j,\rmd x_j) \lambda_{-i} (\rmd t_{-i}) \\
& = \int_{T_{-i}} \int_{X_{-i}} w_i(x_i,x_{-i}, t_i,t_{-i})  \prod_{j \neq i} E^{\lambda_j}\left(f_j|\cG_j \right)(t_j,\rmd x_j) \lambda_{-i} (\rmd t_{-i})  \\
& = \int_{T_{-i}}  w_i(x_i,f_{-i}(t_{-i}), t_i,t_{-i})  \lambda_{-i} (\rmd t_{-i})  \\
& = V_i^{f}(x_i,t_i).
\end{align*}
The third and fifth equalities are due to Equation~(\ref{equa-4}), and the fourth equality holds since $E^{\lambda_i}(g_i|\cG_i) = E^{\lambda_i}(f_i|\cG_i)$ for each $i\in I$.
Thus, $W_i^{(h_i,g)}(t_i) = W_i^{(h_i,f)}(t_i)$ for any $h_i$ and $t_i \in T_i$.

We have
\begin{align*}
U_i(g)
& = \int_{T_i} W_i^{(g_i,g)}(t_i) \lambda_i (\rmd t_i)   = \int_{T_i} W_i^{(g_i,f)}(t_i) \lambda_i (\rmd t_i) \\
& = \int_{T_i} E^{\lambda_i}( W_i^{(g_i,f)}|\cG_i) \lambda_i (\rmd t_i) = \int_{T_i} E^{\lambda_i}(W_i^{(f_i,f)}|\cG_i) \lambda_i (\rmd t_i) \\
& = \int_{T_i} W_i^{(f_i,f)}(t_i) \lambda_i (\rmd t_i) = U_i(f),
\end{align*}
and
\begin{align*}
U_i(h_i,g_{-i})
& = \int_{T_i} W_i^{(h_i,g)}(t_i) \lambda_i (\rmd t_i)   = \int_{T_i} W_i^{(h_i,f)}(t_i) \lambda_i (\rmd t_i) \\
& = U_i(h_i,f_{-i}).
\end{align*}
Thus, $f$ and $g$ are strongly payoff equivalent.

Finally, since $E^{\lambda_i}(c_i^{g_i}|\cG_i) = E^{\lambda_i}(c_i^{f_i}|\cG_i)$, we have
$$ \int_{T_i} c_i^{g_i}(t_i) \lambda_i(\rmd t_i) = \int_{T_i} c_i^{f_i}(t_i) \lambda_i(\rmd t_i) = \int_{T_i}\int_{X_i} c_i(t_i,x_i) f_i(t_i,\rmd x_i) \lambda_i(\rmd t_i) = 1
$$
which implies that $c(t_i,g_i(t_i)) = c_i^{g_i}(t_i) = 1$ for $\lambda_i$-almost all $t_i\in T_i$. That is, $g_i(t_i)\in \supp f_i(t_i)$ for $\lambda_i$-almost all $t_i\in T_i$, $f_i$ and $g_i$ are belief consistent. Since $i$ is arbitrarily chosen, $f$ and $g$ are belief consistent.

Therefore, $g$ is a strong purification of $f$.
\end{proof}

{\small
\singlespacing

\end{document}